\documentclass[10pt, a4paper]{amsart} 
\usepackage[left=2.50cm, right=2.50cm, top=2.50cm, bottom=2.50cm]{geometry} 
\usepackage{amscd, amsfonts, amsmath, amssymb, amsthm, arydshln, fixmath, graphicx, mathrsfs, overpic, tikz, tikz-cd}
\usepackage{hyperref}
\usepackage[all]{xy} 
\usepackage{tikz}
\usepackage{pgfplots}
\usetikzlibrary{arrows}
\usepackage{dsfont}
\makeindex

\theoremstyle{definition} 
\newtheorem{Unity}{Unity}[section] 
\newtheorem*{Definition*}{Definition} 
\newtheorem{Definition}[Unity]{Definition} 

\theoremstyle{plain} 
\newtheorem*{Theorem*}{Theorem}
\newtheorem{Theorem}[Unity]{Theorem}
\newtheorem{Proposition}[Unity]{Proposition}
\newtheorem{Corollary}[Unity]{Corollary}
\newtheorem{Lemma}[Unity]{Lemma}

\theoremstyle{remark} 
\newtheorem*{Remark*}{Remark}
\newtheorem{Remark}[Unity]{Remark}
\newtheorem{Example}[Unity]{Example}

\numberwithin{Unity}{section}

\newcommand{\rk}{\mathrm{rk\,}}

\newcommand{\Rep}{\mathrm{Rep}}

\newcommand{\Spec}{\mathrm{Spec\,}}

\newcommand{\Vect}{\mathfrak{Vect}}

\begin{document}

\title{Base Change Of Local Fundamental Group Schemes}
\author{Lingguang Li}
\address{School of Mathematical Sciences,
Key Laboratory of Intelligent Computing and Applications (Tongji University), Ministry of Education, Shanghai 200092, CHINA}
\email{LiLg@tongji.edu.cn}
\author{Niantao Tian}
\address{School of Mathematical Sciences,
Key Laboratory of Intelligent Computing and Applications (Tongji University), Ministry of Education, Shanghai 200092, CHINA}
\email{tianniantao@tongji.edu.cn}
\begin{abstract}
Let $k$ be a field of characteristic $p>0$, $K/k$ a field extension, $X$ a geometrically reduced connected scheme proper over $k$, and $x\in X(k)$. We study when the natural homomorphism $\pi^{Loc}(X_K,x_K)\rightarrow\pi^{Loc}(X,x)_K$ is an isomorphism. If $\mathcal C^{Loc}(X)$ is saturated, then it is an isomorphism for every algebraic extension $K/k$. If $\pi^{Loc}(X,x)$ is finite, then it is an isomorphism if $K/k$ is purely transcendental or $k$ is algebraically closed. Combining these results, if $\mathcal C^{Loc}(X)$ is saturated and $\pi^{Loc}(X,x)$ is finite, then the natural homomorphism $\pi^{Loc}(X_K,x_K)\rightarrow \pi^{Loc}(X,x)_K$ is an isomorphism. As an application, if $X_{\bar k}$ is smooth and rationally chain connected, then the natural homomorphism $\pi^{Loc}(X_K,x_K)\longrightarrow\pi^{Loc}(X,x)_K$ is an isomorphism.
\end{abstract}
\maketitle

\section{Introduction}
Let $k$ be a field of characteristic $p>0$, $X$ a geometrically reduced connected scheme proper over $k$, and $x\in X(k)$. A vector bundle $E$ on $X$ is called Frobenius trivial if $F_X^{n*}E$ is trivial for some $n\geq1$, where $F_X$ denotes the absolute Frobenius morphism on $X$. These vector bundles form a neutral Tannakian category $\mathcal C^{Loc}(X)$ with the fibre functor at $x$. Its Tannaka group scheme is the local fundamental group scheme $\pi^{Loc}(X,x)$, studied by Mehta and Subramanian \cite{MeSu02,MeSu08}. For a field extension $K/k$, pullback induces a natural homomorphism $\pi^{Loc}(X_K,x_K)\longrightarrow\pi^{Loc}(X,x)_K$. The base change problem asks when this homomorphism is an isomorphism.

Mehta and Subramanian \cite{MeSu02} constructed a counterexample to Nori's base change conjecture on a projective curve with a cuspidal singularity, which is also a counterexample for the local fundamental group scheme. Pauly \cite{Pau07} constructed smooth counterexamples from ordinary curves of genus $2$ in characteristic $2$, and St\"abler \cite{Sta14} constructed smooth projective counterexamples in every positive characteristic.

For a smooth projective variety over an algebraically closed field, Mehta and Subramanian \cite{MeSu08} gave necessary and sufficient conditions for the natural homomorphism $\pi^{Loc}(X_K,x_K)\longrightarrow\pi^{Loc}(X,x)_K$ to be an isomorphism for every extension of algebraically closed fields. Under the same geometric assumptions, Biswas and dos Santos \cite[Section~4.2]{BiDS12} showed that reducedness of the universal Nori torsor is a sufficient condition for this base change property.

In \cite{LiTi26}, we established general Tannakian criteria for base change and applied them to several fundamental group schemes. In particular, the local fundamental group scheme satisfies base change under separable algebraic extensions. The present paper proves base-change isomorphisms for purely inseparable and purely transcendental extensions and gives conditions ensuring base change under arbitrary field extensions. Our base-change results can be summarized as follows; see Propositions~\ref{localgebraicbasechange}, \ref{locfinitepuretranscendentalbasechange}, \ref{locfinitealgclosedbasechange}, and~\ref{rcc-local-basechange}.

\begin{Theorem}
    Let $k$ be a field of characteristic $p>0$, $K/k$ a field extension,  $X$ a geometrically reduced connected scheme proper over $k$, and $x_K\in X_K(K)$ lying over $x\in X(k)$. Then the natural homomorphism $\pi^{Loc}(X_K,x_K)\longrightarrow\pi^{Loc}(X,x)_K$ is an isomorphism if one of the following conditions holds:
    \begin{enumerate}
        \item $K/k$ is purely inseparable and $\mathcal C^{Loc}(X)$ is saturated.
        \item $K/k$ is purely transcendental and $\pi^{Loc}(X,x)$ is finite.
        \item $k$ is algebraically closed and $\pi^{Loc}(X,x)$ is finite.
        \item $\mathcal C^{Loc}(X)$ is saturated and $\pi^{Loc}(X,x)$ is finite.
    \end{enumerate}
\end{Theorem}

\begin{Theorem}[Corollary~\ref{chainconnected-local-basechange}]
    Let $k$ be a field of characteristic $p>0$, $K/k$ a field extension, $X$ a scheme proper over $k$ such that $X_{\bar{k}}$ is smooth and rationally chain connected, and $x\in X(k)$. Then the natural homomorphism $\pi^{Loc}(X_K,x_K) \rightarrow\pi^{Loc}(X,x)_K$ is an isomorphism.
\end{Theorem}

The paper is organized as follows. Section~2 establishes the saturation criterion, its invariance under field extensions, and examples. Section~3 proves the base change results under algebraic and purely transcendental extensions, treats the case over an algebraically closed ground field, and applies the results to rationally chain connected varieties.

\section{Preliminaries}

Let $k$ be a field, $K/k$ a field extension, $X$ a scheme over $k$, $G$ an affine group scheme over $k$ and $x \in X(k)$. We write $\Vect(X)$ for the category of vector bundles on $X$. Consider the Cartesian diagrams
\[
    \begin{aligned}
        \begin{tikzcd}
            X\times_{\Spec k}\Spec K\arrow[r,"\rho"]\arrow[d]& X\arrow[d]\\
            \Spec K\arrow[r]&\Spec k
        \end{tikzcd}&&\begin{tikzcd}
            G\times_{\Spec k}\Spec K\arrow[r]\arrow[d]& G\arrow[d]\\
            \Spec K\arrow[r]&\Spec k.
        \end{tikzcd}
    \end{aligned}
\]
We denote $X_K := X \times_{\Spec k} \Spec K$ and $G_K := G \times_{\Spec k}\Spec K$ for the base changes, and $\rho: X_K \to X$ is the canonical projection and $x_K \in X_K(K)$ is a $K$-valued point lying over $x\in X(k)$. For a coherent sheaf $E$ on $X$, we write $E_K:=E \otimes_k K:=\rho^*E$.

\begin{Definition}
Let $k$ be a field of characteristic $p>0$, $X$ a geometrically reduced connected scheme proper over $k$, $F_X:X\rightarrow X$ the absolute Frobenius morphism, $x\in X(k)$, $E$ a vector bundle on $X$ of rank $r$. We say that $E$ is \textit{Frobenius trivial} if there exists a positive integer $n$ such that $F_X^{n*} E\cong \mathcal{O}_X^{\oplus r}$. We denote by $\mathcal C^{Loc}(X)$ the full subcategory of $\operatorname{Vect}(X)$ consisting of Frobenius-trivial bundles. Together with fibre functor $|_x$, this is a neutral Tannakian category. Its Tannaka group scheme, denoted by $\pi^{Loc}(X,x)$, is called the \textit{local fundamental group scheme}. If $\mathcal{C}^{Loc}(X)$ is closed under extensions of vector bundles on $X$, then we say $\mathcal{C}^{Loc}(X)$ is \textit{saturated}.
\end{Definition}

\begin{Proposition}\label{saturatedLoccriterion}
    Let $k$ be a field of characteristic $p>0$, $X$ a geometrically reduced connected scheme proper over $k$, $F_X:X\rightarrow X$ the absolute Frobenius morphism, and $x_{\bar k}\in X_{\bar k}(\bar k)$ lying over $x\in X(k)$. Then the following conditions are equivalent:
    \begin{enumerate}
        \item $\mathcal C^{Loc}(X)$ is saturated.
        \item The Frobenius action $F_X^*:H^1(X,\mathcal O_X) \rightarrow H^1(X,\mathcal O_X)$ is nilpotent.
        \item $H^1_{\acute{e}t}(X_{\bar k},\mathbb F_p)=0$.
        \item The maximal pro-$p$ quotient of $\pi^{\acute{e}t}_1(X_{\bar k},x_{\bar k})$ is trivial.
        \item Every finite continuous abelian quotient of $\pi^{\acute{e}t}_1(X_{\bar k},x_{\bar k})$ has order prime to $p$.
    \end{enumerate}
\end{Proposition}

\begin{proof}
    Note that $H^0(X,\mathcal O_X)=k$. Indeed, properness implies that
    $H^0(X,\mathcal O_X)$ is a finite-dimensional $k$-algebra. Since $X$ is geometrically reduced, $H^0(X,\mathcal O_X)\otimes_k\bar k\cong
    H^0(X_{\bar k},\mathcal O_{X_{\bar k}})$ is reduced, thus $H^0(X,\mathcal O_X)$ is finite
    \'{e}tale over $k$. Since $X$ is connected, this finite \'{e}tale algebra has no nontrivial idempotents and is therefore a field. The rational point $x$ gives a $k$-algebra homomorphism $H^0(X,\mathcal O_X)\to k$, so $H^0(X,\mathcal O_X)=k$.

    $(1)\Rightarrow(2)$
    If $H^1(X,\mathcal O_X)=0$, there is nothing to prove. Otherwise, choose a $k$-basis $\xi_1,\ldots,\xi_d$ of $H^1(X,\mathcal O_X)$. For each $i$, choose an extension representing $\xi_i$:
    \[
        0\rightarrow\mathcal O_X\rightarrow E_i
        \rightarrow\mathcal O_X\rightarrow0.
    \]
    Since $\mathcal O_X\in\mathcal C^{Loc}(X)$ and
    $\mathcal C^{Loc}(X)$ is saturated, we have
    $E_i\in\mathcal C^{Loc}(X)$.
    There exists an integer $n_i\geq1$ such that $F_X^{n_i*}E_i\cong\mathcal O_X^{\oplus2}$. Pulling back the extension by $F_X^{n_i}$ gives
    an exact sequence
    \[
        0\rightarrow\mathcal O_X
        \rightarrow\mathcal O_X^{\oplus2}
        \rightarrow\mathcal O_X\rightarrow0.
    \]
    Since $H^0(X,\mathcal O_X)=k$, the last map is given by a pair of constants, at least one of which is nonzero. It therefore admits a section, so the sequence splits. Hence $F_X^{n_i*}(\xi_i)=0$.

    Let $m=\max\{n_1,\ldots,n_d\}$. Then for any $a_1,\ldots,a_d\in k$ we have
    \[
        F_X^{m*}\left(\sum_{i=1}^d a_i\xi_i\right)=\sum_{i=1}^d a_i^{p^m}F_X^{m*}(\xi_i)=0.
    \]
    Therefore $F_X^{m*}=0$ on $H^1(X,\mathcal O_X)$.

    $(2)\Rightarrow(1)$ Since the Frobenius action $F_X^*:H^1(X,\mathcal O_X) \rightarrow H^1(X,\mathcal O_X)$ is nilpotent, there exists an integer $m\geq1$ such that $F_X^{m*}=0$ on $H^1(X,\mathcal O_X)$. Consider the exact sequence of vector bundles
    \[
        0\rightarrow E_1\rightarrow E
        \rightarrow E_2\rightarrow0,
    \]
    where $E_1,E_2\in\mathcal C^{Loc}(X)$. There exists an integer $n\geq1$ such that $F_X^{n*}E_i\cong\mathcal O_X^{\oplus r_i}$, $r_i=\rk(E_i)$, $i=1,2$. Applying the pullback by $F_X^n$, we obtain
    \[
        0\rightarrow\mathcal O_X^{\oplus r_1}\rightarrow F_X^{n*}E\rightarrow\mathcal O_X^{\oplus r_2}\rightarrow0.
    \]
    Its extension class belongs to $\operatorname{Ext}^1_{\mathcal O_X}(\mathcal O_X^{\oplus r_2},\mathcal O_X^{\oplus r_1})\cong H^1(X,\mathcal O_X)^{\oplus r_1r_2}$. Applying the pullback by $F_X^m$, we obtain a zero extension class. Consequently, $F_X^{(n+m)*}E\cong\mathcal O_X^{\oplus(r_1+r_2)}$. Thus $E\in\mathcal C^{Loc}(X)$, which implies that $\mathcal C^{Loc}(X)$ is saturated.

    $(2)\Leftrightarrow(3)$ By \cite[Theorem~7.2.11]{Fu11}, $H^1_{\acute{e}t}(X_{\bar k},\mathbb F_p) \cong \ker(F_{X_{\bar k}}^*-\operatorname{id}:H^1(X_{\bar k},\mathcal{O}_{X_{\bar k}})\rightarrow H^1(X_{\bar k},\mathcal{O}_{X_{\bar k}}))$. Let $W=\bigcap_{n\geq0}F_{X_{\bar k}}^{n*}\bigl(H^1(X_{\bar k},\mathcal O_{X_{\bar k}})\bigr)$. Since $\bar k$ is perfect, each image $F_{X_{\bar k}}^{n*}\bigl(H^1(X_{\bar k},\mathcal O_{X_{\bar k}})\bigr)$ is a finite-dimensional $\bar k$-subspace, and hence $W=F_{X_{\bar k}}^{n_0*}\bigl(H^1(X_{\bar k},\mathcal O_{X_{\bar k}})\bigr)$ for some $n_0$ and $F_{X_{\bar k}}^*|_W$ is a $p$-semilinear isomorphism. Note that for any $v\in \ker(F_{X_{\bar k}}^*-\operatorname{id})$, we have $F_{X_{\bar k}}^{n_0*}(v)=v$ and hence $v\in W$, which implies that $\ker(F_{X_{\bar k}}^*-\operatorname{id})=\ker(F_{X_{\bar k}}^*|_W-\operatorname{id}_W)$. Since $F_{X_{\bar k}}^*|_W$ is bijective, \cite[Lemma~4.1.1]{MoRe26} gives a natural isomorphism $\ker\bigl(F_{X_{\bar k}}^*|_W-\operatorname{id}_W\bigr) \otimes_{\mathbb F_p}\bar k \cong W$, which implies that $\dim_{\mathbb F_p} H^1_{\acute{e}t}(X_{\bar k},\mathbb F_p)=\dim_{\bar k}W$. Hence $H^1_{\acute{e}t}(X_{\bar k},\mathbb F_p)=0$ if and only if $W=0$, which is equivalent to $F_{X_{\bar k}}^{n_0*}=0$. Since $F_{X_{\bar k}}^{n*}(\xi\otimes a)=F_X^{n*}(\xi)\otimes a^{p^n}$ for all $\xi\in H^1(X,\mathcal O_X)$, $a\in {\bar k}$, and any $n\geq1$, we have $F_{X}^{n*}=0 \Leftrightarrow F_{X_{\bar k}}^{n*}=0$. Hence $H^1_{\acute{e}t}(X_{\bar k},\mathbb F_p)=0$ if and only if $F_{X}^{n*}=0$ for some $n\geq 1$.

    $(3)\Leftrightarrow(4)\Leftrightarrow(5)$ By \cite[Example~11.3]{Mil13}, $H^1_{\acute{e}t}(X_{\bar k},\mathbb F_p)\cong\operatorname{Hom}_{\mathrm{cont}}(\pi^{\acute{e}t}_1(X_{\bar k},x_{\bar k}),\mathbb Z/p\mathbb Z)$. Hence $(3)$ is equivalent to saying that $\pi^{\acute{e}t}_1(X_{\bar k},x_{\bar k})$ has no continuous quotient group of order $p$. Every nontrivial pro-$p$ group has a nontrivial finite $p$-group quotient, and every nontrivial finite $p$-group has a quotient of order $p$, which implies that $(3)\Leftrightarrow(4)$. A finite abelian group has order divisible by $p$ if and only if it admits a quotient of order $p$. Thus $\pi^{\acute{e}t}_1(X_{\bar k},x_{\bar k})$ has a finite abelian quotient of order divisible by $p$ if and only if it has a quotient of order $p$, which implies that $(3)\Leftrightarrow(5)$.
\end{proof}

\begin{Proposition}\label{saturatedLocbasechange}
    Let $k$ be a field of characteristic $p>0$, $K/k$ a field
    extension, $X$ a geometrically reduced connected scheme
    proper over $k$, and $x_K\in X_K(K)$ lying over $x\in X(k)$.
    Then $\mathcal C^{Loc}(X)$ is saturated if and only if
    $\mathcal C^{Loc}(X_K)$ is saturated.
\end{Proposition}

\begin{proof}
    Since $H^1(X_K,\mathcal O_{X_K})\simeq H^1(X,\mathcal O_X)\otimes_k K$ and $F_{X_K}^{n*}(\xi\otimes a)=F_X^{n*}(\xi)\otimes a^{p^n}$ for all $\xi\in H^1(X,\mathcal O_X)$, $a\in K$, and any $n\geq1$. Hence, for each fixed
    $n\geq 1$, $F_X^{n*}=0 \Leftrightarrow F_{X_K}^{n*}=0$. Then $F_X^*$ is nilpotent on $H^1(X,\mathcal O_X)$ if and only if $F_{X_K}^*$ is nilpotent on $H^1(X_K,\mathcal O_{X_K})$ by Proposition~\ref{saturatedLoccriterion}.
\end{proof}

\begin{Remark}\label{Locsaturatedexamples}
    Let $k$ be a field of characteristic $p>0$, $X$ a geometrically reduced connected scheme proper over $k$, $x\in X(k)$, and $F_X:X\rightarrow X$ the absolute Frobenius morphism. The following criteria and examples hold.
    \begin{enumerate}
        \item If $X$ is a smooth curve, then $\mathcal{C}^{Loc}(X)$ is saturated iff the $p$-rank of the Jacobian of $X$ is zero.
        \item If $X$ is an abelian variety, then $\mathcal{C}^{Loc}(X)$ is saturated iff the $p$-rank of $X$ is zero. In particular, this holds for supersingular abelian varieties.
        \item If $X_{\bar k}$ is smooth and rationally chain connected, then $\mathcal{C}^{Loc}(X)$ is saturated.
    \end{enumerate}
\end{Remark}

\begin{proof}
For an abelian variety $A$ over $k$, write $f(A)=\dim_{\mathbb F_p}A[p](\bar{k})$ for its $p$-rank.

(1) By \cite[\S 3.3, (3.6)--(3.7)]{Ot25}, we have $f(\operatorname{Jac}(X))=\dim_{\mathbb F_p}H^1_{\acute{e}t}(X_{\bar k},\mathbb F_p)$. It follows that $\mathcal C^{Loc}(X)$ is saturated if and only if $f(\operatorname{Jac}(X))=0$ by Proposition~\ref{saturatedLoccriterion}.

(2) By \cite[Exp.~XI, Th\'eor\`eme~2.1]{SGA1}, the maximal pro-$p$ quotient of $\pi_1^{\acute{e}t}(X_{\bar k},0)$ is isomorphic to $\mathbb Z_p^{\oplus f(X)}$. Then it is trivial if and only if $f(X)=0$ and hence $\mathcal C^{Loc}(X)$ is saturated if and only if $f(X)=0$ by Proposition~\ref{saturatedLoccriterion}.

(3) By \cite[Proposition~8.4]{ChLo04}, $\pi^{\acute{e}t}_1(X_{\bar k},x_{\bar k})$ is a finite group of order prime to $p$. Hence $\mathcal{C}^{Loc}(X)$ is saturated by Proposition~\ref{saturatedLoccriterion}.
\end{proof}

\begin{Example}\label{locfinitedoesnotimplysaturated}
    Let $k$ be an algebraically closed field of characteristic $p>0$, $X$ an irreducible projective rational nodal curve of arithmetic genus $1$, $x\in X(k)$, and $F_X:X\rightarrow X$ the absolute Frobenius morphism. Then $\pi^{Loc}(X,x)=1$ but $\mathcal C^{Loc}(X)$ is not saturated.
\end{Example}

\begin{proof}
    Let $z$ be the unique node of $X$, $\nu:\mathbb P_k^1\rightarrow X$ the normalization of $X$. Choose coordinates on $\mathbb P_k^1$ such that $\nu^{-1}(z)=\{0,\infty\}$.

    For any $\mathcal E\in\mathcal C^{Loc}(X)$, there exists $n\geq1$ such that $F_X^{n*}\mathcal E\cong\mathcal O_X^{\oplus r}$ where $\rk(\mathcal{E})=r$. Then we have $F_{\mathbb P_k^1}^{n*}\nu^*\mathcal E\cong\mathcal O_{\mathbb P_k^1}^{\oplus r}$. By \cite[Lemma~6, Theorem~7]{Bu12}, $\nu^*\mathcal E\cong\mathcal O_{\mathbb P_k^1}^{\oplus r}$. Consider the following Cartesian diagram
    \[
        \begin{tikzcd}
            \{0,\infty\}\arrow[r,"\tilde{\eta}"]\arrow[d,"\tilde{\nu}"]& \mathbb P_k^1\arrow[d,"\nu"]\\
            \{z\}\arrow[r,"\eta"]&X.
        \end{tikzcd}
    \]
    By \cite[Theorem~2]{Bu12}, $\mathcal E$ corresponds to the triple $\bigl(\nu^*\mathcal E,\eta^*\mathcal E,m_{\mathcal E}\bigr)$, where $m_{\mathcal E}:\tilde\nu^*\eta^*\mathcal E\xrightarrow{\sim}\tilde{\eta}^*\nu^*\mathcal E$ is the canonical isomorphism. Its restrictions to the two components of $\{0,\infty\}$ are $m_0:\mathcal E|_z\xrightarrow{\sim}\nu^*\mathcal E|_0$ and $m_\infty:\mathcal E|_z\xrightarrow{\sim}\nu^*\mathcal E|_\infty$. Choose a basis of $\mathcal E|_z$, then $m_{\mathcal E}$ is represented by a pair of matrices $M(0)$, $M(\infty)\in\operatorname{GL}_r(k)$. 

By \cite[Lemma~5, Corollary~10]{Bu12}, $F_X^{n*}\mathcal E$ corresponds to the triple $\left( \nu^*F_X^{n*}\mathcal E, \eta^*F_X^{n*}\mathcal E, m_{F_X^{n*}\mathcal E}\right)$ and $m_{F_X^{n*}\mathcal E}$ is represented by a pair of matrices $M'(0):=M(0)^{(p^n)}$, $M'(\infty):=M(\infty)^{(p^n)}$. Fix a basis of $\mathcal O_X^{\oplus r}|_z$, then $\mathcal O_X^{\oplus r}$ corresponds to the triple $\left(\mathcal O_{\mathbb P_k^1}^{\oplus r},\mathcal O_{\{z\}}^{\oplus r},m_{\mathcal O_{X}^{\oplus r}}\right)$ and $m_{\mathcal O_{X}^{\oplus r}}$ is represented by $M''(0)=M''(\infty)=I_r$. Since $F_X^{n*}\mathcal E\cong\mathcal O_X^{\oplus r}$, we may denote this isomorphism by $f:F_X^{n*}\mathcal E\xrightarrow{\cong}\mathcal O_X^{\oplus r}$, then these two triples are isomorphic by \cite[Theorem~2]{Bu12} and we have the following commutative diagram
\[
    \begin{tikzcd}
        \tilde \nu^*\eta^*F_X^{n*}\mathcal E\arrow[r,"m_{F_X^{n*}\mathcal E}"]\arrow[d,"\tilde \nu^*\eta^*(f)"]&\tilde \eta ^*\nu^*F_X^{n*}\mathcal E\arrow[d,"\tilde \eta^*\nu^*(f)"]\\
        \tilde \nu^*\mathcal O_{\{z\}}^{\oplus r}\arrow[r,"m_{\mathcal O_{X}^{\oplus r}}"]&\tilde \eta ^*\mathcal O_{\mathbb P_k^1}^{\oplus r}
    \end{tikzcd}
\]
Since $\nu^*(f)$ is an isomorphism between two trivial vector bundles on $\mathbb P_k^1$, it corresponds to a matrix $P\in \operatorname{GL}_r(k)$. We denote the corresponding matrix of $\eta^*(f)$ by $Q$. Restricting the above commutative diagram to $0$ and $\infty$ respectively, we obtain $P M'(0)=M''(0)Q=Q$, $P M'(\infty)=M''(\infty)Q=Q$. Since $P$ is invertible, we obtain $M(0)^{(p^n)}=M'(0)=M'(\infty)=M(\infty)^{(p^n)}$ and hence $M(0)=M(\infty)$.

Note that the isomorphism $\nu^*\mathcal E\cong\mathcal O_{\mathbb P_k^1}^{\oplus r}$ corresponds to a matrix $P'\in \operatorname{GL}_r(k)$. Then we have the following commutative diagrams
\[
    \begin{tikzcd}
        \tilde \nu^*\eta^*\mathcal E|_0\arrow[r,"M(0)"]\arrow[d,"P'M(0)"]&\tilde \eta ^*\nu^*\mathcal E|_0\arrow[d,"P'"]\\
        \tilde \nu^*\mathcal O_{\{z\}}^{\oplus r}|_0\arrow[r,"I_r"]&\tilde \eta ^*\mathcal O_{\mathbb P_k^1}^{\oplus r}|_0
    \end{tikzcd}
    \qquad
    \begin{tikzcd}
        \tilde \nu^*\eta^*\mathcal E|_\infty\arrow[r,"M(\infty)"]\arrow[d,"P'M(0)"]&\tilde \eta ^*\nu^*\mathcal E|_\infty\arrow[d,"P'"]\\
        \tilde \nu^*\mathcal O_{\{z\}}^{\oplus r}|_\infty\arrow[r,"I_r"]&\tilde \eta ^*\mathcal O_{\mathbb P_k^1}^{\oplus r}|_\infty
    \end{tikzcd}
\]
Since both $P'$ and $P'M(0)$ are invertible, we obtain $\bigl(\nu^*\mathcal E,\eta^*\mathcal E,m_{\mathcal E}\bigr)\cong \left(\mathcal O_{\mathbb P_k^1}^{\oplus r},\mathcal O_{\{z\}}^{\oplus r},m_{\mathcal O_{X}^{\oplus r}}\right)$. Then by \cite[Theorem~2]{Bu12}, we obtain $\mathcal E\cong \mathcal O_{X}^{\oplus r}$ and hence $\pi^{Loc}(X,x)=1$.

    Finally, we have the normalization exact sequence
    \[
        0\rightarrow\mathcal O_X\rightarrow\nu_*\mathcal O_{\mathbb P_k^1}\rightarrow k_z\rightarrow0,
    \]
    where $k_z$ is the skyscraper sheaf at $z$. Then we have the following long exact sequence
    $$\begin{aligned}
        0\rightarrow H^0(X,\mathcal{O}_X)\rightarrow H^0(X,\nu_*\mathcal O_{\mathbb P_k^1})\rightarrow H^0(X,k_z)\xrightarrow{\delta} H^1(X,\mathcal{O}_X)\rightarrow H^1(X,\nu_*\mathcal O_{\mathbb P_k^1})\rightarrow\cdots.
    \end{aligned}$$
    Since $H^0(X,\nu_*\mathcal O_{\mathbb P_k^1})=H^0(\mathbb P_k^1,\mathcal O_{\mathbb P_k^1})$, $H^1(X,\nu_*\mathcal O_{\mathbb P_k^1})=H^1(\mathbb P_k^1,\mathcal O_{\mathbb P_k^1})$, we have long exact sequence
    $$\begin{aligned}
        0\rightarrow k\xrightarrow{\cong} k\xrightarrow{0} k\xrightarrow{\delta} H^1(X,\mathcal{O}_X)\rightarrow 0.
    \end{aligned}$$
    It follows that $\delta$ is an isomorphism and for any $a\in k$, we have $F_X^*(\delta(a))=\delta(a^p)$. Then for any $n\geq 1$, $F_X^{n*}(\delta(1))=F_X^{(n-1)*}(F_X^*(\delta(1)))=F_X^{(n-1)*}(\delta(1))=\delta(1)$, which implies that $F_X^*$ is not nilpotent, and hence $\mathcal C^{Loc}(X)$ is not saturated by Proposition~\ref{saturatedLoccriterion}.
\end{proof}

\section{Base change of local fundamental group schemes}

\begin{Lemma}\label{repeatedlemma}
    Let $k$ be a field, $P$ a scheme proper over $k$ with $H^0(P,\mathcal O_P)=k$, $U$ an integral scheme of finite type over $k$, $f:P\times_k U\rightarrow U$ the natural projection, $\mathcal F$ a vector bundle of rank $r$ on $P\times_k U$. Then there exists a nonempty open subscheme $U'\subseteq U$, $f_*\mathcal F$ is a free sheaf on $U'$ and the natural map $(f_*\mathcal F)|_s \rightarrow H^0(P_{\kappa(s)},{\mathcal F}|_{P_{\kappa(s)}})$ is an isomorphism for any $s\in U'$. If $\mathcal F|_{P_{\kappa(s)}}\cong \mathcal O_{{P_{\kappa(s)}}}^{\oplus r}$ for some closed point $s\in U'$, then $\mathcal F|_{P_{\kappa(\eta)}}\cong \mathcal O_{{P_{\kappa(\eta)}}}^{\oplus r}$, where $\eta$ is the generic point of $U$.
\end{Lemma}

\begin{proof}
    Since $f$ is proper and ${\mathcal F}$ is flat over $U$. By generic freeness and generic cohomology-and-base-change imply, we can shrink $U$ such that $f_*{\mathcal F}$ is free on $U$. Let $\eta\in U$ be the generic point. Consider the Cartesian diagrams 
    \[
        \begin{tikzcd}
            P_{\kappa(\eta)}\arrow[r]\arrow[d] &P\times_k U\arrow[d,"f"]\\
            \Spec(\kappa(\eta))\arrow[r]&U,
        \end{tikzcd}
    \]
    by cohomology and flat base change, we have $(f_*\mathcal F)|_\eta\cong H^0(P_{\kappa(\eta)},\mathcal F|_{P_{\kappa(\eta)}})$. Then there exists nonempty open subscheme $U'\subseteq U$ such that for any $s\in U'$, the natural map $(f_*\mathcal F)|_s \rightarrow H^0(P_{\kappa(s)},{\mathcal F}|_{P_{\kappa(s)}})$ is an isomorphism.

    Since $H^0(P,\mathcal O_P)=k$, we have $H^0(P_{\kappa(u)},\mathcal O_{P_{\kappa(u)}})=\kappa(u)$ for any $u\in U'$ by field-extension base change. Suppose $\mathcal F|_{P_{\kappa(s)}}\cong \mathcal O_{{P_{\kappa(s)}}}^{\oplus r}$ for a closed point $s\in U'$. Then $(f_*\mathcal F)|_{s}\cong H^0(P_{\kappa(s)},{\mathcal F}|_{P_{\kappa(s)}})\cong \kappa(s)^{\oplus r}$, so $f_*\mathcal F$ is free of rank $r$. Consider the adjunction evaluation map $\alpha:f^*f_*{\mathcal F}\rightarrow{\mathcal F}$ on $P\times_k U$. Since $\alpha$ is a morphism between vector bundles of the same rank, the non-isomorphism locus of $\alpha$ is $\operatorname{Supp}(\operatorname{coker}\alpha)$. Restricting $\alpha$ to $P_{\kappa(s)}$, we obtain a natural isomorphism
    $$f^*f_*{\mathcal F}|_{P_{\kappa(s)}}\cong H^0(P_{\kappa(s)},{\mathcal F}|_{P_{\kappa(s)}})\otimes_{{\kappa(s)}}\mathcal O_{P_{\kappa(s)}}\xrightarrow{\cong}{\mathcal F}|_{P_{\kappa(s)}}\cong\mathcal O_{P_{\kappa(s)}}^{\oplus r}.$$
    It follows that $\alpha|_{y}$ is an isomorphism for any $y\in P_{\kappa(s)}$. Hence $P_{\kappa(s)}\bigcap \operatorname{Supp}(\operatorname{coker}\alpha)= \varnothing$. Since $f$ is proper, $f(\operatorname{Supp}(\operatorname{coker}\alpha))$ is closed in $U$ and $s\in U\setminus f(\operatorname{Supp}(\operatorname{coker}\alpha))$. Then $\eta\in U\setminus f(\operatorname{Supp}(\operatorname{coker}\alpha))$ and restricting $\alpha$ to $P_{\kappa(\eta)}$ gives $\mathcal O_{P_{\kappa(\eta)}}^{\oplus r}\cong \mathcal F|_{P_{\kappa(\eta)}}$. 
\end{proof}

\subsection{Base change under algebraic extensions}

\begin{Proposition}\label{localgebraicbasechange}
    Let $k$ be a field of characteristic $p>0$, $K/k$ an algebraic field extension, $X$ a geometrically reduced connected scheme proper over $k$, and $x_K\in X_K(K)$ lying over $x\in X(k)$. Suppose one of the following conditions is satisfied:
    \begin{enumerate}
        \item $K/k$ is separable.
        \item $K/k$ is purely inseparable and $\mathcal{C}^{Loc}(X)$ is saturated.
    \end{enumerate}
    Then the natural homomorphism $\pi^{Loc}(X_K,x_K)\rightarrow\pi^{Loc}(X,x)_K$ is an isomorphism.
\end{Proposition}

\begin{proof}
    (1) It follows by \cite[Theorem~4.11]{LiTi26}.

    (2) Suppose $K/k$ is finite purely inseparable. Since $\mathcal{C}^{Loc}(X_K)$ is saturated by Proposition~\ref{saturatedLocbasechange}, the homomorphism $\pi^{Loc}(X_K,x_K)\rightarrow\pi^{Loc}(X,x)_K$ is an isomorphism by \cite[Proposition~3.12, Proposition~4.7]{LiTi26}. 
    
    Suppose $K/k$ is infinite purely inseparable. Let $\mathcal{E}\in \mathcal{C}^{Loc}(X_K)$, there exists a finite extension $L/k$ contained in $K$ and $E_L\in\Vect(X_L)$ such that $\mathcal{E}\cong E_L\otimes_L K$. By \cite[Proposition~4.7]{LiTi26}, we have $E_L\in\mathcal{C}^{Loc}(X_L)$. Since $L/k$ is finite, we have $\pi^{Loc}(X_L,x_L)\cong \pi^{Loc}(X,x)_L$, then
    $$E_L\in\mathcal{C}^{Loc}(X_L)\simeq\Rep_L^f(\pi^{Loc}(X,x)_L).$$
    By \cite[Theorem~3.3]{LiTi26}, there exists $F\in\mathcal{C}^{Loc}(X)$ such that $E_L\hookrightarrow F\otimes_k L$ and thus 
    $$\mathcal{E}\cong E_L\otimes_L K\hookrightarrow (F\otimes_k L)\otimes_L K\cong F\otimes_k K.$$
    Hence the homomorphism $\pi^{Loc}(X_K,x_K)\rightarrow \pi^{Loc}(X,x)_K$ is an isomorphism by \cite[Theorem~3.3]{LiTi26}.
\end{proof}

\begin{Example} \label{counterexample}
The curve and the failure of base change below are from \cite[Example~4.13]{LiTi26}. We determine the two local fundamental group schemes. Let $k=\mathbb F_2(t)$, $K=k(s)$, $s^2=t$, $C$ be the smooth projective plane cubic over $k$ with affine equation $y^2+xy=x^3+t$, $o$ its point at infinity and $F_{C_K}:C_K\rightarrow C_K$ the absolute Frobenius morphism. Then $K/k$ is purely inseparable, $\pi^{Loc}(C,o)=1$, $\pi^{Loc}(C_K,o_K)\cong\mu_{2,K}$ and hence the natural homomorphism $\pi^{Loc}(C_K,o_K) \rightarrow\pi^{Loc}(C,o)_K$ is not an isomorphism. 
\end{Example}

\begin{proof}
    The Weierstrass invariants are $\Delta=t$ and $c_4=1$, so $j(C)=t^{-1}$. Fix $L\in\{k,K\}$ and let $G$ be any finite quotient of $\pi^{Loc}(C_L,o_L)$, which is also a finite local quotient of $\pi^N(C_L,o_L)$. By \cite[Remark~3, p.~264]{Nor83}, we have $\pi^N(C_L,o_L)\cong\varprojlim_{m\ge1}C_L[m]$. Then there is a faithfully flat homomorphism $C_L[m]\to G$ for some $m$. Let $\ker$ be its kernel. Multiplication by $m$ on $C_L$ induces an isogeny $C_L/\ker\rightarrow C_L$ of elliptic curves with kernel $C_L[m]/\ker\cong G$. As $G$ is local, this isogeny is purely inseparable. If its degree is $2^a$, it is the $a$-fold relative Frobenius of $C_L/\ker$ up to an $L$-isomorphism of its target \href{https://stacks.math.columbia.edu/tag/0CCZ}{Tag 0CCZ}. In particular, $|G|=2^a$, $j(C_L)=j(C_L/\ker)^{2^a}$.

    Over $k=\mathbb F_2(t)$, applying the discrete $t$-adic valuation to $j(C)=j(C/\ker)^{2^a}$ gives $-1=2^a v_t(j(C/\ker))$, so $a=0$. Thus all finite quotients of $\pi^{Loc}(C,o)$ are trivial, and $\pi^{Loc}(C,o)=1$.

    Over $K=\mathbb F_2(s)$, we have $j(C_K)=s^{-2}$. Applying the discrete $s$-adic valuation to $j(C_K)=j(C_K/\ker)^{2^a}$ gives $-2=2^a v_s(j(C_K/\ker))$, so $a\le1$. Hence every finite quotient of $\pi^{Loc}(C_K,o_K)$ has order at most $2$. Let $P=(0,s)\in C_K(K)$, $\mathcal{L}=\mathcal O_{C_K}([P]-[o_K])$. The inverse of a point $(x,y)$ on $C$ is $(x,y+x)$. Hence $P=-P$ and $P\ne o_K$. It follows that $\mathcal{L}^{\otimes2}\cong\mathcal O_{C_K}$, $F_{C_K}^*\mathcal{L}\cong\mathcal O_{C_K}$. Thus $\mathcal L\in\mathcal C^{Loc}(C_K)$. Then the Tannakian subcategory generated by $\mathcal L$ corresponds to a finite quotient of $\pi^{Loc}(C_K,o_K)$ which is nontrivial and hence of order $2$. The finite local group scheme of order $2$ is either $\mu_{2,K}$ or $\alpha_{2,K}$, up to isomorphism, and the finite dimensional $K$-representation category of $\alpha_{2,K}$ has a unique irreducible object which is $K$. It follows that the finite quotient of $\pi^{Loc}(C_K,o_K)$ corresponds to the Tannakian subcategory generated by $\mathcal L$ is isomorphic to $\mu_{2,K}$. Let $G'$ be any finite quotient of $\pi^{\mathrm{Loc}}(C_K,o_K)$. There exists a common finite quotient $H'$ admitting faithfully flat morphisms onto both $G'$ and $\mu_{2,K}$. Since $|H'|\le 2$, the quotient $H'\twoheadrightarrow\mu_{2,K}$ implies that $H'\simeq\mu_{2,K}$. Since $\pi^{\mathrm{Loc}}(C_K,o_K)$ is profinite, it follows that $\pi^{\mathrm{Loc}}(C_K,o_K)\simeq\mu_{2,K}$.
\end{proof}

\subsection{Base change under purely transcendental extensions}
\begin{Proposition}\label{locfinitepuretranscendentalbasechange}
    Let $k$ be a field of characteristic $p>0$, $K/k$ a purely transcendental field extension, $X$ a geometrically reduced connected scheme proper over $k$, $F_X:X\rightarrow X$ the absolute Frobenius morphism, and $x_K\in X_K(K)$ lying over $x\in X(k)$. Suppose $\pi^{Loc}(X,x)$ is a finite group scheme. Then the natural homomorphism $\pi^{Loc}(X_K,x_K)\rightarrow\pi^{Loc}(X,x)_K$ is an isomorphism.
\end{Proposition}

\begin{proof}
    Let $q:P\rightarrow X$ be the universal pointed $\pi^{Loc}(X,x)$-torsor. Since $\pi^{Loc}(X,x)$ is a finite group scheme over $k$, the associated bundle corresponding to the regular representation $k[\pi^{Loc}(X,x)]$ is $q_*\mathcal O_P$, so $q_*\mathcal O_P\in\mathcal C^{Loc}(X)$. Moreover, the Tannakian equivalence
    $\operatorname{Rep}_k^f(\pi^{Loc}(X,x))\simeq\mathcal C^{Loc}(X)$ gives $H^0(P,\mathcal O_P)\simeq\operatorname{Hom}_X(\mathcal O_X,q_*\mathcal O_P)\simeq\operatorname{Hom}_G(k,k[\pi^{Loc}(X,x)])\simeq k$.

    For any $\mathcal E\in\mathcal C^{Loc}(X_K)$, there exists an integer $n\geq1$ such that $F_{X_K}^{n*}\mathcal E \cong\mathcal O_{X_K}^{\oplus r}$ where $r=\rk(\mathcal E)$. Suppose first that $K/k$ has finite transcendence degree, and write $K=k(t_1,\ldots,t_d)$, where $t_1,\ldots,t_d$ are algebraically independent over $k$. By spreading out $\mathcal E$ together with this trivialization, there exist a nonempty affine open subset $U\subseteq\mathbb A_k^d$ and a vector bundle $\mathcal E_U$ on $X_U$ whose generic fibre is $\mathcal E$, such that $F_{X_U}^{n*}\mathcal E_U\cong\mathcal O_{X_U}^{\oplus r}$. Then $\mathcal E_U|_{X_{\kappa(s)}}\in\mathcal C^{Loc}(X_{\kappa(s)})$ for every $s\in U$, where $\kappa (s)$ is the residue field of $s$.

    Let $q_U:P_U\rightarrow X_U$ be the base change of $q$, and $f:P_U\rightarrow U$ the projection. According to Lemma~\ref{repeatedlemma}, we can shrink $U$ such that $f_*{q_U^*\mathcal E_U}$ is free and the natural map $(f_*q_U^*\mathcal E_U)|_s \rightarrow H^0(P_{\kappa(s)},{q_U^*\mathcal E_U}|_{P_{\kappa(s)}})$ is an isomorphism for any $s\in U$. Choose a closed point $s\in U$ such that $\kappa(s)/k$ is finite separable. Such a point exists: if $k$ is infinite, then $U$ has a $k$-rational point; if $k$ is finite, then every closed point of $U$ has residue field finite separable over $k$. By Proposition~\ref{localgebraicbasechange}, we have $\pi^{Loc}(X_{\kappa(s)},x_{\kappa(s)})\cong\pi^{Loc}(X,x)_{\kappa(s)}$. Hence $q_{\kappa(s)}:P_{\kappa(s)}\rightarrow X_{\kappa(s)}$ is the universal local torsor. Since $\mathcal E_U|_{X_{\kappa(s)}}\in\mathcal C^{Loc}(X_{\kappa(s)})$, we obtain ${q_U^*\mathcal E_U}|_{P_{\kappa(s)}}=q_{\kappa(s)}^*(\mathcal E_U|_{X_{\kappa(s)}})\cong\mathcal O_{P_{\kappa(s)}}^{\oplus r}$. According to Lemma~\ref{repeatedlemma}, we have ${q_U^*\mathcal E_U}|_{P_{\kappa(\eta)}}=q_{\kappa(\eta)}^*(\mathcal E_U|_{X_{\kappa(\eta)}})\cong\mathcal O_{P_{\kappa(\eta)}}^{\oplus r}$, where $\eta$ is the generic point of $U$, i.e. $\mathcal O_{P_K}^{\oplus r}\cong q_K^*\mathcal E$.

    For an arbitrary purely transcendental extension $K/k$, and any $\mathcal E\in\mathcal C^{Loc}(X_K)$ with $F_{X_K}^{n*}\mathcal E \cong\mathcal O_{X_K}^{\oplus r}$ where $r=\rk(\mathcal E)$, the bundle and its Frobenius trivialization descend to $X_{k(t_1,\ldots,t_d)}$ for finitely many elements $t_1,\ldots,t_d$ of a transcendence basis over $k$. Applying the preceding argument and then extending scalars to $K$, we again obtain $q_K^*\mathcal E\cong\mathcal O_{P_K}^{\oplus r}$. Hence $q_K:P_K\rightarrow X_K$ trivializes every object of $\mathcal C^{Loc}(X_K)$.

    Since $q_K$ is faithfully flat, the adjunction $\mathcal  E\rightarrow(q_K)_*q_K^*\mathcal E$ is injective. Consequently, $\mathcal E \hookrightarrow(q_K)_*q_K^*\mathcal E\cong\bigl((q_K)_*\mathcal O_{P_K}\bigr)^{\oplus r}\cong\bigl(q_*\mathcal O_P\bigr)^{\oplus r}\otimes_k K$. Then by \cite[Theorem~3.3]{LiTi26}, the natural homomorphism $\pi^{Loc}(X_K,x_K)\rightarrow\pi^{Loc}(X,x)_K$ is an isomorphism. 
\end{proof}

\subsection{Base change over algebraically closed fields}

\begin{Proposition}\label{locfinitealgclosedbasechange}
    Let $k$ be an algebraically closed field of characteristic $p>0$, $K/k$ a field extension,  $X$ a geometrically reduced connected scheme proper over $k$, and $x_K\in X_K(K)$ lying over $x\in X(k)$. Suppose that $\pi^{Loc}(X,x)$ is a finite group scheme. Then the natural homomorphism $\pi^{Loc}(X_K,x_K) \rightarrow \pi^{Loc}(X,x)_K$ is an isomorphism.
\end{Proposition}

\begin{proof}
    Let $q:P\rightarrow X$ be the universal pointed $\pi^{Loc}(X,x)$-torsor. Since $\pi^{Loc}(X,x)$ is a finite group scheme over $k$, we have $q_*\mathcal O_P\in\mathcal C^{Loc}(X)$ and $H^0(P,\mathcal O_P)=k$ as explained in the proof of
    Proposition~\ref{locfinitepuretranscendentalbasechange}.

    For any $\mathcal E\in\mathcal C^{Loc}(X_K)$, there exists an integer $n\geq1$ such that $F_{X_K}^{n*}\mathcal E\cong\mathcal O_{X_K}^{\oplus r}$, where $r=\operatorname{rk}(\mathcal E)$. By spreading out $\mathcal E$ together with this trivialization, there exist a finitely generated $k$-subalgebra $R\subseteq K$ and a vector bundle $\mathcal E_U$ on $X_U$, where $U=\operatorname{Spec}R$, such that $\mathcal E_U\otimes_R K\cong\mathcal E$, $F_{X_U}^{n*}\mathcal E_U\cong\mathcal O_{X_U}^{\oplus r}$. Then $\mathcal E_U|_{X_{\kappa(s)}}\in\mathcal C^{Loc}(X_{\kappa(s)})$ for every $s\in U$.

    Let $q_U:P_U\rightarrow X_U$ be the base change of $q$, and $f:P_U\rightarrow U$ the projection. According to Lemma~\ref{repeatedlemma}, we can shrink $U$ such that $f_*{q_U^*\mathcal E_U}$ is free and the natural map $(f_*q_U^*\mathcal E_U)|_s \rightarrow H^0(P_{\kappa(s)},{q_U^*\mathcal E_U}|_{P_{\kappa(s)}})$ is an isomorphism for any $s\in U$. Choose a closed point $s\in U$. Since $k$ is algebraically closed and $U$ is of finite type over $k$, we have $\kappa(s)=k$. Hence $\mathcal E_U|_{X_{\kappa(s)}}\in\mathcal C^{Loc}(X_{\kappa(s)})$ and we have $(q_U^*\mathcal E_U)|_{P_{\kappa(s)}}\cong q_{\kappa(s)}^*(\mathcal E_U|_{X_{\kappa(s)}})\cong\mathcal O_{P_{\kappa(s)}}^{\oplus r}$. According to Lemma~\ref{repeatedlemma}, we have ${q_U^*\mathcal E_U}|_{P_{\kappa(\eta)}}=q_{\kappa(\eta)}^*(\mathcal E_U|_{X_{\kappa(\eta)}})\cong\mathcal O_{P_{\kappa(\eta)}}^{\oplus r}$. After base change to $K$, we obtain $\mathcal O_{P_{K}}^{\oplus r}\cong q_{K}^*\mathcal E$. Hence $q_K:P_K\rightarrow X_K$ trivializes every object of $\mathcal C^{Loc}(X_K)$.

    Finally, note that $q_*\mathcal O_P\in\mathcal C^{Loc}(X)$ because $\pi^{Loc}(X,x)$ is finite. Since $q_K$ is faithfully flat, the adjunction $\mathcal  E\rightarrow(q_K)_*q_K^*\mathcal E$ is injective. Consequently, $\mathcal E \hookrightarrow(q_K)_*q_K^*\mathcal E\cong\bigl((q_K)_*\mathcal O_{P_K}\bigr)^{\oplus r}\cong\bigl(q_*\mathcal O_P\bigr)^{\oplus r}\otimes_k K$. Then by \cite[Theorem~3.3]{LiTi26}, the natural homomorphism $\pi^{Loc}(X_K,x_K)\rightarrow\pi^{Loc}(X,x)_K$ is an isomorphism. 
\end{proof}
\subsection{Applications}

\begin{Proposition}\label{rcc-local-basechange}
Let $k$ be a field of characteristic $p>0$, $K/k$ a field extension, $X$ a geometrically reduced connected scheme proper over $k$, and $x\in X(k)$. Suppose $\mathcal C^{Loc}(X)$ is saturated and $\pi^{Loc}(X,x)$ is a finite group scheme. Then the natural homomorphism $\pi^{Loc}(X_K,x_K) \rightarrow\pi^{Loc}(X,x)_K$ is an isomorphism.
\end{Proposition}

\begin{proof}
    Consider the factorization $K/L_s/L/k$, where $L/k$ is purely transcendental, $K/L$ is algebraic and $L_s$ is the subfield of $K$ consisting of the elements separable algebraic over $L$. Then $L_s/L$ is separable algebraic and $K/L_s$ is purely inseparable algebraic. By Proposition~\ref{locfinitepuretranscendentalbasechange}, the natural homomorphism $\pi^{Loc}(X_{L},x_{L})\rightarrow \pi^{Loc}(X,x)_{L}$ is an isomorphism and $\pi^{Loc}(X_L,x_L)$ remains finite. Proposition~\ref{saturatedLocbasechange} shows that
    $\mathcal C^{\mathrm{Loc}}(X_{L_s})$ is saturated. Then by Proposition~\ref{localgebraicbasechange}, the natural homomorphisms $\pi^{Loc}(X_{L_s},x_{L_s})\rightarrow \pi^{Loc}(X_{L},x_{L})_{L_s}$ and $\pi^{Loc}(X_K,x_K)\rightarrow \pi^{Loc}(X_{L_s},x_{L_s})_K$ are isomorphisms. Hence the homomorphism $\pi^{Loc}(X_K,x_K) \rightarrow\pi^{Loc}(X,x)_K$ is an isomorphism.
\end{proof}

\begin{Corollary}\label{chainconnected-local-basechange}
    Let $k$ be a field of characteristic $p>0$, $K/k$ a field extension, $X$ a scheme proper over $k$ such that $X_{\bar{k}}$ is smooth and rationally chain connected, and $x\in X(k)$. Then the natural homomorphism $\pi^{Loc}(X_K,x_K) \rightarrow\pi^{Loc}(X,x)_K$ is an isomorphism.
\end{Corollary}

\begin{proof}
    Since $X_{\bar k}$ is smooth and rationally chain connected proper over $\bar{k}$, by \cite[Proposition~8.4]{ChLo04}, $\pi^{\acute{e}t}_1(X_{\bar k},x_{\bar k})$ is a finite group of order prime to $p$. Then by Proposition~\ref{saturatedLoccriterion} and Proposition~\ref{saturatedLocbasechange}, $\mathcal C^{Loc}(X)$ is saturated, and hence the natural homomorphism $\pi^{Loc}(X_{\bar k},x_{\bar k})\rightarrow\pi^{Loc}(X,x)_{\bar k}$ is an isomorphism under factorization ${\bar k}/k^{sep}/k$ by Proposition~\ref{localgebraicbasechange}. Since $X_{\bar k}$ is normal and rationally chain connected, \cite[Theorem~3.10]{AnBi16} shows that $\pi^{Loc}(X_{\bar k},x_{\bar k})$ is a finite group scheme over $\bar k$. Then the group scheme $\pi^{Loc}(X,x)$ is also a finite group scheme over $k$ because finiteness descends under the faithfully flat field extension ${\bar k}/k$. According to Proposition~\ref{rcc-local-basechange}, the natural homomorphism $\pi^{Loc}(X_K,x_K) \rightarrow\pi^{Loc}(X,x)_K$ is an isomorphism.
\end{proof}

\noindent\textbf{AI disclosure.}
The authors used ChatGPT to assist in constructing the examples in Example~\ref{locfinitedoesnotimplysaturated} and Example~\ref{counterexample}. Apart from this author-directed assistance, all ideas and arguments are due to the authors. The authors independently verified all arguments and take full responsibility for the paper.

\end{document}